\documentclass[a4paper,12pt]{amsart}
\usepackage{amssymb}
\usepackage{ifthen}
\usepackage{graphicx}

\newtheorem{thm}{Theorem}[section]
\newtheorem{cor}{Corollary}[section]
\newtheorem{lem}{Lemma}[section]

\newtheorem{rem}{Remark}[section]
\newtheorem{example}{Example}[section]

\newtheorem{prob}{Problem}[section]

\newtheorem{conj}{Conjecture}
\theoremstyle{definition}

\numberwithin{equation}{section}

\newcounter {own}
\def\theown {\thesection       .\arabic{own}}

\newenvironment{pf}[1][]{%
 \vskip 3mm
 \noindent
 \ifthenelse{\equal{#1}{}}%
  {{\slshape Proof. }}%
  {{\slshape #1.} }%
 }%
{\qed\bigskip}

\newcounter{alphabet}
\newcounter{tmp}
\newenvironment{Thm}[1][]{\refstepcounter{alphabet}%
\bigskip%
\noindent%
{\bf Theorem \Alph{alphabet}}%
\ifthenelse{\equal{#1}{}}{}{ (#1)}%
{\bf .} \itshape}{\vskip 8pt}

\makeatletter
\newcommand{\reff}[1]{\@ifundefined{r@#1}{}{\setcounter{tmp}{\ref{#1}}\Alph{tmp}}}
\makeatother

\newcommand{\IC}{{\mathbb C}}

\newcommand{\ID}{{\mathbb D}}

\newcommand{\D}{{\mathbb D}}

\def\be{\begin{equation}}
\def\ee{\end{equation}}

\newcommand{\bee}{\begin{enumerate}}
\newcommand{\eee}{\end{enumerate}}

\newcommand{\blem}{\begin{lem}}
\newcommand{\elem}{\end{lem}}
\newcommand{\bthm}{\begin{thm}}
\newcommand{\ethm}{\end{thm}}
\newcommand{\bcor}{\begin{cor}}
\newcommand{\ecor}{\end{cor}}
\newcommand{\beg}{\begin{example}}
\newcommand{\eeg}{\end{example}}
\newcommand{\begs}{\begin{examples}}
\newcommand{\eegs}{\end{examples}}
\newcommand{\bdefe}{\begin{defin}}
\newcommand{\edefe}{\end{defin}}
\newcommand{\bprob}{\begin{prob}}
\newcommand{\eprob}{\end{prob}}
\newcommand{\bei}{\begin{itemize}}
\newcommand{\eei}{\end{itemize}}

\newcommand{\bcon}{\begin{conj}}
\newcommand{\econ}{\end{conj}}
\newcommand{\bcons}{\begin{conjs}}
\newcommand{\econs}{\end{conjs}}
\newcommand{\bprop}{\begin{propo}}
\newcommand{\eprop}{\end{propo}}
\newcommand{\br}{\begin{rem}}
\newcommand{\er}{\end{rem}}
\newcommand{\brs}{\begin{rems}}
\newcommand{\ers}{\end{rems}}
\newcommand{\bo}{\begin{obser}}
\newcommand{\eo}{\end{obser}}
\newcommand{\bos}{\begin{obsers}}
\newcommand{\eos}{\end{obsers}}
\newcommand{\bpf}{\begin{pf}}
\newcommand{\epf}{\end{pf}}
\newcommand{\ba}{\begin{array}}
\newcommand{\ea}{\end{array}}
\newcommand{\beq}{\begin{eqnarray}}
\newcommand{\beqq}{\begin{eqnarray*}}
\newcommand{\eeq}{\end{eqnarray}}
\newcommand{\eeqq}{\end{eqnarray*}}

\newcounter{minutes}
\divide\time by 60
\newcounter{hours}
\multiply\time by 60 \addtocounter{minutes}{-\time}

\begin{document}
\bibliographystyle{amsplain}
\title{Lipschitz property of harmonic mappings with respect to the pseudo-hyperbolic metric}
\thanks{
\indent
File:~\jobname .tex,
printed: \number\day-\number\month-\number\year,
\thehours.\ifnum\theminutes<10{0}\fi\theminutes
}
\author{Jie Huang}
\address{Jie Huang, School of Mathematical Sciences, Huaqiao
University, Quanzhou 362021, China; Mathematics with Computer Science Program, Guangdong Technion, 241 Daxue Road, Jinping District, Shantou, Guangdong 515063, People's Republic of China and Department of Mathematics, Technion -- Israel Institute of Technology, Haifa 32000, Israel.
}
\email{jie.huang@gtiit.edu.cn}

\author{Antti Rasila${}^*$}

\address{Antti Rasila, Mathematics with Computer Science Program, Guangdong Technion, 241 Daxue Road, Jinping District, Shantou, Guangdong 515063, People's Republic of China and Department of Mathematics, Technion -- Israel Institute of Technology, Haifa 32000, Israel.}
\email{antti.rasila@iki.fi; antti.rasila@gtiit.edu.cn}

\author{Jian-Feng Zhu}
\address{Jian-Feng Zhu, School of Mathematical Sciences, Huaqiao
University, Quanzhou 362021, China and Department of Mathematics,
Shantou University, Shantou, Guangdong 515063, China.}
\email{flandy@hqu.edu.cn}


\begin{abstract}
 In this paper, we show that harmonic Bloch mappings are Lipschitz continuous with respect to the pseudo-hyperbolic metric. This result improves the corresponding result of $($\cite[Theorem 1]{Proceeding-AMS}$)$. Furthermore, we prove the similar property for harmonic quasiregular Bloch-type mappings. 
\end{abstract}

\thanks{
The research of the authors were supported by NNSF of China (No. 11501220, 11971182, 11971124), NSF of Fujian Province (No. 2016J01020, 2019J0101), NSF of Guangdong province (No. 2021A1515010326), the Promotion Program for Young and Middle-aged Teachers in Science, and Technology Research of Huaqiao University (ZQN-PY402). * Corresponding author.
}

\maketitle \pagestyle{myheadings}
\markboth{J. Huang, A. Rasila, and J.-F. Zhu}{Lipschitz property of harmonic mappings with respect to pseudo-hyperbolic metric}

\section{Introduction}
Let $\ID=\{z\in\mathbb{C}:\; |z|<1\}$ be the unit disk, $\mathbb{T}=\{z\in\mathbb{C}:\; |z|=1\}$ the unit circle, and $\overline{\ID}$ the closure of $\mathbb{D}$, i.e., $\overline{\ID}=\ID\cup \mathbb{T}$.
For $z\in \mathbb{C}$, the partial derivatives of a complex-valued function $f$ are defined by:
\beqq\label{eq1.1}
f_z=\frac{1}{2}\left(f_x-if_y\right)\;\;\mbox{and}\;\;f_{\overline{z}}=\frac{1}{2}\left(f_x+if_y\right).
\eeqq
For $z=re^{i\theta}\in \mathbb{C}$ and $\alpha\in[0,2\pi]$, the directional derivative of $f$ is defined by
\beqq\label{eq1.2}
\partial_{\alpha}f(z)=\lim\limits_{r\rightarrow0^+}\frac{f(z+re^{i\alpha})-f(z)}{r}=e^{i\alpha}f_z(z)+e^{-i\alpha}f_{\bar{z}}(z).
\eeqq
Then
\beqq\label{eq1.3}
\Lambda_f(z):=\max\limits_{0\leq\alpha\leq2\pi}|\partial_{\alpha}f(z)|=|f_z(z)|+|f_{\bar{z}}(z)|
\eeqq
and
\beqq\label{eq1.4}
\lambda_f(z):=\min\limits_{0\leq\alpha\leq2\pi}|\partial_{\alpha}f(z)|=\big||f_z(z)|-|f_{\bar{z}}(z)|\big|.
\eeqq

A complex-valued function $f(z)$ of the class $C^2$ is said to be harmonic mapping, if it satisfies
$\Delta f=4f_{z\bar{z}}=0$. Moreover, it was shown in \cite{Lewy1936} that a function $f$ is locally univalent and
sense-preserving in $\ID$ if and only if its Jacobian $J_f(z)=|f_z(z)|^2-|f_{\bar{z}}(z)|^2>0$, i.e.,
the dilatation of $f$
$$|\omega_f(z)|=\frac{|f_{\bar{z}}(z)|}{|f_z(z)|}<1$$
in $\ID$.
Assume that $f(z)$ is a harmonic mapping defined in a simply connected domain $\Omega\subseteq \IC$.
Then $f(z)$ has the canonical decomposition $f(z)=h(z)+\overline{g(z)}$, where $h(z)$ and $g(z)$ are analytic in $\Omega$.
For a sense-preserving harmonic mapping $f(z)$ in $\ID$, let
$$\omega (z)=\frac{g'(z)}{h'(z)}$$
be the (second) complex dilatation of $f$. Then $\omega (z)$ is holomorphic mapping of $\ID$ and
$$\|\omega (z)\|_{\infty }:= \sup_{z\in \ID}|\omega (z)|\leq 1.$$
In this paper, we consider locally univalent and sense-preserving harmonic mappings in $\ID$.

It is worth of noting that the composition $f\circ \varphi$ of a harmonic mapping $f$ with a conformal mapping $\varphi$ is
a harmonic mapping. Therefore, it is sufficient to consider harmonic mappings defined in the unit disk $\ID$.
However, $\varphi\circ f$ is not harmonic in general.

Fix $w\in\ID$ and let $\varphi_w$ be the M\"{o}bius transformation of $\ID$, that is,
$$\varphi_w(z)=\frac{w-z}{1-\bar{w}z}, \ \ \ \mbox{where}\ \ \ z\in\ID.$$
The pseudo-hyperbolic distance on $\ID$ is defined by
$$
\rho(z, w)=|\varphi_w(z)|.
$$
The pseudo-hyperbolic distance is M\"{o}bius invariant, that is,
$$\rho(g(z), g(w))=\rho(z, w),$$
for all $g\in$ Aut$(\ID)$, the M\"{o}bius automorphisms of $\ID$. It has the following useful property (which will be used in proving Theorem \ref{HZ-2-thm-1} below):
\be\label{hz-2018-eq-2}
1-\rho(z, w)^2=\frac{(1-|z|^2)(1-|w|^2)}{|1-\bar{z}w|^2}=(1-|w|^2)|\varphi_z'(w)|.
\ee
\noindent{\bf Definition 1.1.} {\it We call a function $h$ a { Bloch function} (and write $h\in B$) $h$ is analytic in $\ID$ and }
$$\|h\|_{B}=\sup\limits_{z\in\ID}(1-|z|^2)|h'(z)|<\infty.$$
The above formula defines a seminorm, and the Bloch functions form a complex Banach space $\mathcal{B}$ with the norm
$$\|h\|_{\mathcal{B}}=|h(0)|+\|h\|_{B}.$$

A mapping $f(z)$ is said to be Lipschitz (resp. co-Lipschitz) in $\ID$ if there exists a constant $L$ such that
the following inequality
$$
\frac{|z_1-z_2|}{L}\leq|f(z_1)-f(z_2)|\ \ \ \mbox{(resp.}\ \ \ |f(z_1)-f(z_2)|\leq L|z_1-z_2|)
$$
holds for all $z_1$, $z_2\in\ID$, where $L\geq1$ is called the Lipschitz constant. The function $f$ is said to be bi-Lipschitz if $f$ is Lipschitz and co-Lipschitz.

It is easy to show that the condition $h\in B$ does not ensure that $h$ is a Lipschitz mapping. For example, take
$h(z)=\log(1-z^2)$, for $z\in\ID$. Then $h\in B$ since $\|h\|_B=(1-|z|^2)|h'(z)|\leq2$.
However, for $z_1=x\in (0, 1)$, choose arbitrary small $t>0$ such that $z_2=x+t\in (0, 1)$. Then
$$\left|\frac{h(z_1)-h(z_2)}{z_1-z_2}\right|\geq\frac{1}{1-x}\rightarrow\infty,\ \ \  \mbox{as}\ \ \ x\rightarrow1.$$
This shows that $h$ is not a Lipschitz mapping.

Let
$$H(z, w):=C_{\varphi_w} h(z)=h\circ\varphi_w(z),$$
where $z, w\in\ID$. Then $|H'(0, w)|=(1-|w|^2)|h'(w)|$.

Lipschitz and Bloch type spaces of harmonic mappings and their generalizations have been by several authors in, for example, \cite{CLR20,CPR1,CPR2,Chenshaolin-11,Colonna,Samy-RIM}.
In \cite{Proceeding-AMS}, P. Ghatage, J. Yan, and D. Zheng showed that $(1-|w|^2)|h'(w)|$ is a Lipschitz function with respect to pseudo-hyperbolic metric. They also used this result to study the composition operators $C_{\varphi_w}$ on the Bloch space.
In fact, they proved the following Theorem.

\begin{Thm}\label{Thm-1}$($\cite[Theorem 1]{Proceeding-AMS}$)$
Let $h$ be in the Bloch space. Then the following inequality
$$\big|(1-|z|^2)|h'(z)|-(1-|w|^2)|h'(w)|\big|\leq3.31\rho(z, w)\|h\|_{B},$$
holds for all $z, w\in\ID$.
\end{Thm}
In this paper, our primary goal is to improve the above Theorem \reff{Thm-1} in the case of harmonic Bloch-type mappings.

\noindent{\bf Definition 1.2.} {\it A harmonic mapping $f$ in $\ID$ is called a harmonic Bloch mapping (denote by $f\in B_h$) if }
$$\|f\|_{B_h}=\sup\limits_{z\in\ID}(1-|z|^2)\Lambda_f(z)<\infty.$$

This defines a seminorm, and the space equipped with the norm
$$\|f\|_{\mathcal{B}_h}=|f(0)|+\|f\|_{B_h}$$
is called the {\it harmonic Bloch space}, denote by $\mathcal{B}_h$. It is a Banach space.
Clearly, $f=h+\bar{g}\in B_h$ if and only if $h\in B$ and $g\in B$, since
$$\max\{\|h\|_B, \|g\|_{B}\}\leq \|f\|_{B_h}\leq \|h\|_B+\|g\|_B.$$
The harmonic Bloch space was studied by Colonna in \cite{Colonna} as a generalization
of the classical Bloch space. We refer to \cite{Anderson,Bonk,Chenshaolin-11,Pommerenke-70,Pommerenke} and
the references therein for more information on $B_h$ and $\mathcal{B}_h$.

Motivated by a number of well-known results on
analytic Bloch functions, in \cite{Efraimids,Samy-RIM} the authors introduced the {\it harmonic Bloch-type mappings}, which was defined as follows:

\noindent{\bf Definition 1.3.} {\it A harmonic mapping $f$ in $\ID$ is called a harmonic Bloch-type mapping if}
$$\|f\|_{B_h^*}=\sup\limits_{z\in\ID}(1-|z|^2)\sqrt{|J_f(z)|}<\infty.$$

We denote this class of functions by $B_h^*$ and call the quality
$$\|f\|_{\mathcal{B}_h^*}=|f(0)|+\|f\|_{B_h^*}$$
the Bloch-type pseudo-norm of $f$.
It is easy to see that $B_h\subseteq B_h^*$, because $\sqrt{|J_f(z)|}\leq \Lambda_f(z)$, for each $z\in\ID$.

In this paper, we improve Theorem \reff{Thm-1} as follows:

\begin{thm}\label{HZ-2-thm-1}
Let $f$ be in $B_h$ space. Then the following inequality
$$
\big|(1-|z|^2)\Lambda_f(z)-(1-|w|^2)\Lambda_f(w)\big|\leq5.7174\rho(z, w)\|f\|_{B_h},
$$
holds for all $z, w\in\ID$.
\end{thm}

Suppose $f(z)$ is a sense-preserving harmonic mapping
of $\ID$ into a domain $\Omega\subseteq\mathbb{C}$. Then $f(z)$ is
a harmonic $K$-quasiregular mapping, if
$$
K(f):=\sup\limits_{z\in\ID}\frac{|f_z(z)|+|f_{\bar{z}}(z)|}{|f_z(z)|-|f_{\bar{z}}(z)|}\leq K,
$$
where $K\geq1$ is a constant.

We show in Lemma \ref{HZ-2-lemma2} below that if $f$ is a harmonic quasiregular mapping in $\ID$, then $f\in B_h^*$ if and only if $f\in B_h$.
Moreover, by using quasiregularity, we generalize Theorem \reff{Thm-1} as follows.

\begin{thm}\label{HZ-2-thm-2}
Let $f$ be a harmonic K-quasiregular mapping in $\ID$ and in $B_h^*$. Then the following inequality
\be\label{HZ-2-thm-2-eq-1}
\bigg|(1-|z|^2)\sqrt{J_f(z)}-(1-|w|^2)\sqrt{J_f(w)}\bigg|\leq3.6920(K+1)\rho(z, w)\|f\|_{B_h^*},
\ee
holds for all $z, w\in\ID$.
\end{thm}

\section{Auxiliary results}
In this section, we will prove three lemmas. Lemma \ref{HZ-2-lemma1} and Lemma \ref{HZ-2-lemma3} will be used in proving Theorem \ref{HZ-2-thm-1},
and Lemma \ref{HZ-2-lemma2} will be used in proving Theorem \ref{HZ-2-thm-2}.
We start with the following.

\begin{lem}\label{HZ-2-lemma1}
Let $f$ be in $B_h$ or in $B_h^*$. The respective pseudo-norms of $f$ are M\"{o}bius invariant.
\end{lem}
\bpf
For $z, w\in\ID$, let $$\lambda=\varphi_w(z)=\frac{w-z}{1-\overline{w}z}.$$ 
Then $$z=\varphi_w(\lambda)=\frac{w-\lambda}{1-\overline{w}\lambda}.$$
Elementary calculation leads to
$$(1-|z|^2)|\varphi_w'(z)|=1-|\lambda|^2.$$
Thus
\begin{eqnarray*}
\|f\circ\varphi_w\|_{B_h}&= &\sup\limits_{z\in\ID}(1-|z|^2)\big(|f_\lambda(\lambda)||\varphi_w'(z)|+|f_{\bar{\lambda}}(\lambda)||\varphi_w'(z)|\big)\\ \nonumber
&=&\sup\limits_{z\in\ID}(1-|\lambda|^2)\Lambda_f(\lambda)=\|f\|_{B_h}.
\end{eqnarray*}
Similarly, we have $\|f\circ\varphi_w\|_{B_h^*}=\|f\|_{B_h^*}.$
This completes the proof of Lemma \ref{HZ-2-lemma1}.
\epf

The following lemma shows that if $f$ is a harmonic quasiregular mapping of $\ID$, then $f\in B_h^*$ if and only if
$f\in B_h$.

\begin{lem}\label{HZ-2-lemma2}
Let $f$ be a harmonic quasiregular mapping in $\ID$. Then $f\in B_h$ if and only if $f\in B_h^*$. Moreover, we have
$$
\|f\|_{B_h^*}\leq\|f\|_{B_h}\leq\sqrt{K}\|f\|_{B_h^*}.
$$
\end{lem}
\bpf
Suppose $f\in B_h$. Since $\sqrt{J_f(z)} \leq \Lambda _f(z)$, it is easy to see that
\begin{eqnarray*}
\sup_{z\in \D}(1-|z|^2)\sqrt{J_f(z)}\leq \sup_{z\in \D}(1-|z|^2)\Lambda _f(z)=\|f\|_{B_h}.
\end{eqnarray*}
This implies that $f\in B_h^*$ and $ \|f\|_{B_h^*}\leq\|f\|_{B_h}$.

On the other hand, suppose $f\in B_h^*$. The assumption that $f$ is a harmonic quasiregular mapping of $\ID$, ensures that $f$
has the canonical decomposition $f=h+\bar{g}$, where $h$ and $g$ are analytic in $\ID$ , and
$$\frac{|h'|+|g'|}{|h'|-|g'|}=\frac{\Lambda _f^2}{J_f}\leq K.$$
This implies that $\Lambda _f\leq \sqrt K \sqrt {J_f}$, and thus, $f\in B_h$. Moreover, we have
$$\|f\|_{B_h}\leq \sqrt K\|f\|_{B^*_h}.$$
The proof of Lemma \ref{HZ-2-lemma2} is complete.
\epf

\begin{lem}\label{HZ-2-lemma3}
Let $h$ be in the Bloch space $B$.
For $z, w\in\ID$, let $\zeta =\varphi_w(z)$ and $g=C_{\varphi_w} h=h\circ\varphi_w$. If $|\zeta|\leq\frac{1}{3}$,
then
$$
(1-|\zeta|^2)\left|g'(\zeta)-g'(0)\right|\leq c_1|\zeta|\|h\|_B,
$$
where $c_1\approx 2.6920$ the minimal value of 
\[
\psi(r)=\frac{1+r^2/9}{r(1-r^2)}, \qquad 0<r<1.
\]
\end{lem}
\bpf
For $z, w\in\ID$, recall that
\[
\zeta=\varphi_w(z)=\frac{w-z}{1-\bar{w}z}.
\]
Then
$$(1-|w|^2)h'(w)=(h\circ\varphi_w)'(0).$$
Let $g=h\circ\varphi_w$. Then we may rewrite the above equation as
$$g'(0)=(1-|w|^2)h'(w).$$
Note that for any $w\in\ID$, by using the mean-value property for analytic functions, we obtain 
$$
(1-|w|^2)|g''(w)|=|(g'\circ \varphi _w)'(0)|=\frac{1}{2\pi r}\left|\int_0^{2\pi }g'\circ \varphi _w(re^{i\theta })e^{-i\theta }d\theta \right|,
$$
where $0<r<1$.
By estimating the integral and noting that $\|g\|_B=\|h\|_B$, one has
\begin{eqnarray*}
&&\left|\int_0^{2\pi }g'\circ \varphi _w(re^{i\theta })e^{-i\theta }d\theta \right|\leq \|g\|_B\int _0^{2\pi }\frac{1}{1-|\varphi _w(re^{i\theta })|^2}d\theta \\ \nonumber
&=&\|h\|_B\int _0^{2\pi }\frac{|1-\overline{w}re^{i\theta }|^2}{(1-|w|^2)(1-r^2)}d\theta.
\end{eqnarray*}
It follows from the inequality
$$\int_0^{2\pi}\mbox{Re}(\overline{w}re^{i\theta})d\theta=0$$
and the above discussions that the equality
$$(1-|w|^2)^2|g''(w)|\leq \|h\|_B\frac{1+r^2|w|^2}{r(1-r^2)}$$
holds for any $0<r<1$.

Now, consider the function $$\psi(r)=\frac{1+r^2/9}{r(1-r^2)},$$ where $0<r<1$.
Denote by $c_1\approx 2.6920$ the minimal value of $\psi(r)$.
Thus if $|w|\leq \frac{1}{3}$, then
\begin{eqnarray*}
(1-|w|^2)^2|g''(w)|\leq c_1\|h\|_B.
\end{eqnarray*}

By using the inequality $$|g'(\zeta )-g'(0)|\leq \int_0^1|g''(t\zeta )||\lambda |dt$$ and the assumption that $|\zeta|\leq\frac{1}{3}$, we have
\begin{eqnarray*}
|g'(\zeta )-g'(0)|&\leq &c_1\|h\|_B\int _0^1\frac{|\zeta |}{(1-t^2|\zeta|^2)^2}dt\\ \nonumber
&=&c_1\|h\|_B\int _0^{|\zeta|}\frac{ds}{(1-s^2)^2}\\ \nonumber
&=&c_1\|h\|_B\cdot \frac{1}{4}\bigg[\frac{2|\zeta|}{1-|\zeta|^2}+\ln \frac{1+|\zeta |}{1-|\zeta |}\bigg].
\end{eqnarray*}
This implies that
$$(1-|\zeta |^2)|g'(\zeta )-g'(0)|\leq c_1|\zeta |\|h\|_B,$$
because $$(1-|\zeta|^2)\ln \frac{1+|\zeta |}{1-|\zeta |}\leq2|\zeta|,$$
completing the proof of Lemma \ref{HZ-2-lemma3}.
\epf

\section{Proof of main results}
\subsection*{Proof of Theorem \ref{HZ-2-thm-1}}
Let $$\zeta =\varphi_w(z)=\frac{w-z}{1-\overline{w}z},$$ where $z,w\in \ID$, and let $\psi=f\circ \varphi _w$, where $f\in B_h$. Then
\be\label{April-8}(1-|w|^2)\Lambda _f(w)=|\psi_z(0)|+|\psi_{\bar{z}}(0)|=\Lambda _{\psi}(0).\ee
On the other hand,  as $$z=\frac{w-\lambda}{1-\overline{w}\zeta}=\varphi_w(\lambda),$$ it follows from (\ref{hz-2018-eq-2}) that
\begin{eqnarray*}
(1-|z|^2)\Lambda _f(z)&=&(1-|\varphi_w(\zeta )|^2)\Lambda _f(\varphi _w(\zeta ))\\ \nonumber
&=&(1-|\zeta |^2 )|\varphi' _w(\zeta )|\Lambda _f(\varphi _w(\zeta ))\\ \nonumber
&=&(1-|\zeta |^2)\Lambda _{\psi }(\zeta ).
\end{eqnarray*}
Hence,
\be\label{main-1}\left|(1-|z|^2)\Lambda _f(z)-(1-|w|^2)\Lambda _f(w)\right|=\left|(1-|\zeta|^2)\Lambda _{\psi }(\zeta )-\Lambda _{\psi }(0)\right|.\ee
We now divide our proof into two cases.
\begin{description}
  \item[Case 1] Let $|\zeta| \leq \frac{1}{3}$.
\end{description}

First, it follows from (\ref{main-1}) that
$$\left|(1-|z|^2)\Lambda _f(z)-(1-|w|^2)\Lambda _f(w)\right|\leq|\zeta |^2\Lambda _{\psi }(0)+(1-|\zeta|^2)\left|\Lambda _{\psi}(\zeta)-\Lambda _{\psi}(0)\right|.$$
According to (\ref{April-8}), Definition 1.2 and Lemma \ref{HZ-2-lemma1}, we have
\be\label{6-2-1}
\Lambda _{\psi}(0)\leq \|\psi \|_{B_h}=\|f\circ \varphi _w\|_{B_h}=\|f\|_{B_h}. \ee

Next, we estimate $\left|\Lambda _{\psi}(\zeta )-\Lambda _{\psi}(0)\right|$.
According to \cite{Du-04}, we can write $\psi$ as $\psi =h\circ \varphi _w+\overline{g\circ \varphi _w}$.
By letting $H=h\circ \varphi _w$ and $G=g\circ \varphi _w$, both are analytic in $\ID$, we have
$$\left|\Lambda _{\psi}(\zeta )-\Lambda _{\psi}(0)\right|\leq |H'(\zeta )-H'(0)|+|G'(\zeta)-G'(0)|.$$
The assumption $f=h+\bar{g}\in B_h$ ensures that $h$ and $g$ are in $B$.
It follows from Lemma \ref{HZ-2-lemma3} that
\be\label{6-2-2}(1-|\zeta |^2)|H'(\zeta )-H'(0)|\leq c_1 |\zeta |\|h\|_B\leq c_1|\zeta |\|f\|_{B_h}\ee
and
\be\label{6-2-3}(1-|\zeta |^2)|G'(\zeta )-G'(0)|\leq c_1|\zeta |\|g\|_B\leq c_1|\zeta |\|f\|_{B_h}.\ee
Thus, $$(1-|\zeta |^2)\left|\Lambda _{\psi}(\zeta )-\Lambda _{\psi}(0)\right|\leq 2c_1|\zeta|\|f\|_{B_h}.$$
Applying (\ref{6-2-1}), (\ref{6-2-2}) and (\ref{6-2-3}), we have
\begin{eqnarray}\label{case-1}
|\zeta |^2\Lambda _{\psi }(0)+(1-|\zeta |^2)\left|\Lambda _{\psi}(\zeta )-\Lambda _{\psi}(0)\right|
&\leq &|\zeta |^2\|f\|_{B_h}+2c_1|\zeta \|f\|_{B_h}\\ \nonumber
&\leq &c_2|\zeta|\|f\|_{B_h},
\end{eqnarray}
where $c_2\approx 5.7174$.
\begin{description}
  \item[Case 2] Let $\frac{1}{3}<|\zeta|<1$.
\end{description}
In this case, we have $3 |\zeta|>1$. Then
\begin{eqnarray}\label{case-2}
\left|(1-|\zeta |^2)\Lambda _{\psi }(\zeta )-\Lambda _{\psi }(0)\right|&\leq &\max\big\{(1-|\zeta |^2)\Lambda _{\psi }(\zeta ),\Lambda _{\psi }(0)\big\}\\ \nonumber
&\leq &\|\psi \|_{B_h}\\ \nonumber
&< &3|\zeta |\|f\|_{B_h}.
\end{eqnarray}
Combining (\ref{case-1}) and (\ref{case-2}) show that
\be\label{6-2-4}\left|(1-|\zeta|^2)\Lambda _{\psi }(\zeta )-\Lambda _{\psi }(0)\right|\leq c_2|\zeta |\|f\|_{B_h},\ee
holds for all $|\zeta |<1$.

Based on the above discussions, the desired inequality
$$\big|(1-|z|^2)\Lambda_f(z)-(1-|w|^2)\Lambda_f(w)\big|\leq c_2\rho(z, w)\|f\|_{B_h},$$
follows from (\ref{main-1}) and (\ref{6-2-4}), because $\rho (z,w)=|\zeta|$.
This completes the proof of Theorem  \ref{HZ-2-thm-1}.
\qed
\subsection*{Proof of Theorem \ref{HZ-2-thm-2}}
Let $\zeta =\varphi_w(z)$ and $\psi=f\circ \varphi _w$, where $z,w\in \ID$ and $f=h+\bar{g}$ is a harmonic quasiregular mapping in $\ID$,
where $h$ and $g$ are analytic in $\ID$.
Then $J_f(z)=|h'(z)|^2-|g'(z)|^2>0$, and
$$\sqrt{J_{\psi }(0)}=(1-|w|^2)\sqrt{J_f(w)}.$$
Moreover, it follows from (\ref{hz-2018-eq-2}) that
\begin{eqnarray*}
(1-|z|^2)\sqrt{J_f(z)}=(1-|\zeta |^2)\sqrt{J_{\psi }(\zeta )}.
\end{eqnarray*}
Hence,
\be\label{main-2} \big|(1-|z|^2)\sqrt{J_f(z)}-(1-|w|^2)\sqrt{J_f(w)}\big|=\big|(1-|\zeta |^2)\sqrt{J_{\psi }(\zeta )}-\sqrt{J_{\psi }(0)}\big|.\ee
Following the proof of Theorem \ref{HZ-2-thm-1}, we now divide our proof into two cases.
\begin{description}
  \item[Case 1] Let $|\zeta| \leq \frac{1}{3}$.
\end{description}

First, it follows from (\ref{main-2}) that
\begin{eqnarray*}
&&\left|(1-|z|^2)\sqrt{J_f(z)}-(1-|w|^2)\sqrt{J_f(w)}\right|\\
\nonumber
&\leq& |\zeta |^2\sqrt{J_{\psi }(0)}+(1-|\zeta |^2)\left|\sqrt{J_{\psi}(\zeta )}-\sqrt{J_{\psi}(0)}\right|
\end{eqnarray*}
By using Definition 1.3, we have
$$\sqrt{J_{\psi}(0)}\leq \sup\limits_{w\in\ID}(1-|w|^2)\sqrt{J_f(w)}=\|f\|_{B_h^*}.$$

Next, we estimate $\left|\sqrt{J_{\psi}(\zeta )}-\sqrt{J_{\psi}(0)}\right|$ as follows:

By letting $H=h\circ \varphi _w$ and $G=g\circ \varphi _w$, we have
\begin{eqnarray*}
\big|\sqrt{J_{\psi}(\zeta )}-\sqrt{J_{\psi}(0)}\big|&=&\frac{\big|J_{\psi}(\zeta )-J_{\psi}(0)\big|}{\sqrt{J_{\psi}(\zeta )}+\sqrt{J_{\psi}(0)}}\\ \nonumber
&\leq &\frac{\big||H'(\zeta )|+|H'(0)|\big|\cdot\big||H'(\zeta )|-|H'(0)|\big|}{|H'(\zeta )|\sqrt{1-|\omega_\psi(\zeta)|^2}+|H'(0)|\sqrt{1-|\omega_\psi(0)|^2}}\\ \nonumber
&&+\frac{\big||G'(\zeta )|+|G'(0)|\big|\cdot\big||G'(\zeta )|-|G'(0)|\big|}{|H'(\zeta )|\sqrt{1-|\omega_\psi(\zeta)|^2}+|H'(0)|\sqrt{1-|\omega_\psi(0)|^2}},
\end{eqnarray*}
where $\omega_\psi=G'/H'$.
Because $f=h+\bar{g}$ is a harmonic quasiregular mapping of $\ID$,
$$\|\omega_f\|_{\infty}=\sup\limits_{z\in\ID}\frac{|g'(z)|}{|h'(z)|}\leq k,$$
where $k=\frac{K-1}{K+1}<1$.
Direct calculations lead to
$$\|\omega_\psi\|_{\infty}=\sup\limits_{z\in\ID}\frac{|G'(z)|}{|H'(z)|}\leq k.$$
Then
\begin{eqnarray*}\label{6-2-thm2-1}
\frac{\big||H'(\zeta )|+|H'(0)|\big|\cdot\big||H'(\zeta )|-|H'(0)|\big|}{|H'(\zeta )|\sqrt{1-|\omega_\psi(\zeta)|^2}+|H'(0)|\sqrt{1-|\omega_\psi(0)|^2}}\leq
\frac{|H'(\zeta )-H'(0)|}{\sqrt{1-k^2}},
\end{eqnarray*}
and
\begin{eqnarray*}\label{6-2-thm2-2}
\frac{\big||G'(\zeta )|+|G'(0)|\big|\cdot\big||G'(\zeta )|-|G'(0)|\big|}{|H'(\zeta )|\sqrt{1-|\omega_\psi(\zeta)|^2}+|H'(0)|\sqrt{1-|\omega_\psi(0)|^2}}
\leq \frac{|G'(\zeta )-G'(0)|}{\sqrt{1-k^2}}.
\end{eqnarray*}
These show that
\be\label{April-9-1}\left|\sqrt{J_{\psi}(\zeta )}-\sqrt{J_{\psi}(0)}\right|\leq\frac{|H'(\zeta )-H'(0)|+|G'(\zeta )-G'(0)|}{\sqrt{1-k^2}}.\ee

Moreover, because $f\in B_h^*$, we see from Lemma \ref{HZ-2-lemma2} that
$$\|f\|_{B_h}\leq \sqrt{K}\|f\|_{B_h^*}.$$
Therefore, it follows from Lemma \ref{HZ-2-lemma3} that
\begin{eqnarray*}\label{kqc-1}
(1-|\zeta|^2)|H'(\zeta )-H'(0)|&\leq &c_1|\zeta |\|H\|_B\\ \nonumber
&\leq &c_1|\zeta|\|f\|_{B_h}\\ \nonumber
&\leq &c_1|\zeta |\sqrt{K}\|f\|_{B_h^*},
\end{eqnarray*}
and similarly,
$$
(1-|\zeta |^2)|G'(\zeta )-G'(0)|\leq c_1|\zeta |\sqrt{K}\|f\|_{B_h^*}.
$$
Combining the above inequalities and (\ref{April-9-1}) yields
\begin{eqnarray*}
(1-|\zeta |^2)\left|\sqrt{J_{\psi}(\zeta )}-\sqrt{J_{\psi}(0)}\right|
&\leq &\frac {2c_1\sqrt K|\zeta|\|f\|_{B_h^*}}{\sqrt{1-k^2}}\\ \nonumber
&=& c_1(K+1)|\zeta|\|f\|_{B_h^*}.
\end{eqnarray*}
Hence, for $|\zeta|\leq\frac{1}{3}$, one has
\begin{eqnarray}\label{main-thm2-1}
&&|\zeta|^2\sqrt{J_{\psi }(0)}+(1-|\zeta |^2)\big|\sqrt{J_{\psi}(\zeta)}-\sqrt{J_{\psi}(0)} \big|\\ \nonumber
&\leq &c_1(K+1)|\zeta |\|f\|_{B_h^*}+|\zeta |^2\|f\|_{B_h^*}\\ \nonumber
&\leq&c_3(K+1)|\zeta |\|f\|_{B_h^*},
\end{eqnarray}
where $c_3\approx 3.6920$.

\begin{description}
  \item[Case 2] Let $\frac{1}{3}<|\zeta|<1$.
\end{description}
In this case, we have $2\sqrt 2 |\zeta |>1$. Then
\begin{eqnarray}\label{main-thm2-2}
\big|(1-|\zeta |^2)\sqrt{J_{\psi }(\zeta )}-\sqrt{J_{\psi }(0)}\big|&\leq &\max\big\{((1-|\zeta |^2)\sqrt{J_{\psi }(\zeta )},\sqrt{J_{\psi }(0)}\big\}\\ \nonumber
&\leq &\|\psi \|_{B_h^*}\\ \nonumber
&< &3|\zeta |\|f\|_{B_h^*}.
\end{eqnarray}
The desired inequality (\ref{HZ-2-thm-2-eq-1}) follows from (\ref{main-thm2-1}) and (\ref{main-thm2-2}).
This completes the proof of Theorem \ref{HZ-2-thm-2}.
\qed

\end{document}